\documentclass[12pt,reqno]{amsart}                  
\usepackage{amssymb}
\usepackage{mathrsfs}

\newtheorem{theorem}{Theorem}
\newtheorem{lemma}[theorem]{Lemma}
\begin{document}
\title[A proof of the Dalang-Morton-Willinger theorem]{A proof of the Dalang-Morton-Willinger theorem}

\author{Dmitry B. Rokhlin}

\address{D.B. Rokhlin,
Faculty of Mathematics, Mechanics and Computer Sciences,
              Southern Federal University, 
Mil'chakova str., 8a, 344090, Rostov-on-Don, Russia}          
\email{rokhlin@math.rsu.ru}                   
\newcommand{\ri}{{\rm ri\,}}
\newcommand{\lin}{{\rm lin\,}}
\newcommand{\cl}{{\rm cl}}
\newcommand{\epi}{{\rm epi\,}}
\newcommand{\cone}{{\rm cone\,}}
\newcommand{\intern}{{\rm int\,}}
\newcommand{\conv}{{\rm conv\,}}
\newcommand{\dom}{{\rm dom\,}}
\newcommand{\supp}{{\rm supp\,}}
\newcommand{\esssup}{{\rm ess\,sup}}

\begin{abstract}  We give a new proof of the Dalang-Morton-Willinger theorem, relating the no-arbitrage condition in stochastic securities market models to the existence of an equivalent martingale measure with bounded density for a $d$-dimensional stochastic sequence $(S_n)_{n=0}^N$ of stock prices. Roughly speaking, the proof is reduced to the assertion that under the no-arbitrage condition for $N=1$ and $S\in L^1$ there exists a strictly positive linear fucntional on $L^1$, which is bounded from above on a special subset of the subspace $K\subset L^1$ of investor's gains. 
\end{abstract}

\maketitle

\section{Introduction} 
Let $(\Omega,\mathscr F,\mathsf P)$ be a probability space, endowed with a discrete-time filtration $\mathbb F=(\mathscr F_n)_{n=0}^N$, $\mathscr F_N=\mathscr F$. Consider a $d$-dimensional stochastic process $S=(S_n)_{n=0}^N$, adapted to the filtration $\mathbb F$, and a $d$-dimensional $\mathbb F$-predictable process $\gamma=(\gamma_n)_{n=1}^N$, that is, $\gamma_n$ is $\mathscr F_{n-1}$-measurable.
In the customary securities market model $S^i_n$ describes the discounted price of $i$th stock and $\gamma^i_n$ corresponds to the number of stock units in investor's portfolio at time moment $n$. The discounted gain process is given by
\begin{displaymath}
G_n^\gamma=\sum_{k=1}^n(\gamma_k,\Delta S_k),\ \ \ \Delta S_k=S_k-S_{k-1},\ \ n=1,\dots,N,
\end{displaymath}
where $(\cdot,\cdot)$ is the scalar product in $\mathbb R^d$.

Let's recall the classical Dalang-Morton-Willinger theorem \cite{DalMorWil90}, \cite{Shi98} (ch.V, \S 2e).
As usual, we say that the \emph{No Arbitrage (NA)} condition is satisfied if the inequality $G_N^\gamma\ge 0$ a.s. (with respect to the measure $\mathsf P$) implies that $G_N^\gamma=0$ a.s. A probability measure $\mathsf Q$ on $\mathscr F$ is called a \emph{martingale measure} if the process $S$ is a $\mathsf Q$-martingale. Measures $\mathsf P$ and $\mathsf Q$ are called \emph{equivalent} if their null sets are the same. Denote by $\varkappa_{n-1}(\omega)$ the \emph{support of the regular conditional distribution} $\mathsf P_{n-1}(\omega,dx)$ of the random vector $\Delta S_n$ with respect to $\mathscr F_{n-1}$:
\begin{displaymath}
\varkappa_{n-1}(\omega)=\{x\in\mathbb R^d:\mathsf P_{n-1}(\omega,B_\varepsilon(x))>0\ \textrm{for all }\varepsilon>0\},
\end{displaymath}
where $B_\varepsilon(x)\subset\mathbb R^d$ is the ball of center $x$ and radius $\varepsilon$. 
\begin{theorem}[Dalang-Morton-Willinger]
The following conditions are equivalent:
\begin{itemize}
\item[(i)] NA;
\item[(ii)] there exists an equivalent to $\mathsf P$ martingale measure $\mathsf Q$ with a.s. bounded density $z=d\mathsf Q/d\mathsf P$;
\item[(iii)] the relative interior of the convex hull of $\varkappa_{n-1}$ contains the origin a.s., $n=1,\dots,N$.
\end{itemize}
\end{theorem}

At this degree of generality Theorem 1 was proved in \cite{DalMorWil90}. Besides the original one, several alternative proofs \cite{Sch92}, \cite{KabKra94}, \cite{Rog94}, \cite{JacShi98}, \cite{KabStr01} were proposed. The "difficult"
part of this theorem is the implication (i) $\Longrightarrow$ (ii) or (iii) $\Longrightarrow$ (ii), though, in fact, it is sufficient to consider the case of integrable process $S$ and $N=1$.   

In \cite{DalMorWil90} and \cite{JacShi98} the desired martingale measure is constructed at first in the case of trivial $\sigma$-algebra $\mathscr F_0$ and then in general case with the use of subtle measurable selection arguments.  Another approach, based on the closedness in probability (under the NA condition) of the cone $C$ of random variables, dominated by investor's gains $G_N^\gamma$, was proposed in \cite{Sch92}. By the Kreps-Yan theorem this property implies the existence of a strictly positive essentially bounded function $g$ 
such that $\mathsf E(wg)\le 0$, $w\in C\cap L^1$, where $\mathsf E$ is the expectation with respect to $\mathsf P$ and $L^1$ is the space of $\mathsf P$-integrable functions. The function $g$, up to the normalization constant, coincides with the $\mathsf P$-density of an equivalent martingale measure.  The subsequent work in this direction \cite{KabKra94}, \cite{KabStr01} allowed to simplify the proof of the closedness of $C$, as well as the proof of the existence of $g$. We also have mentioned the paper \cite{Rog94} (see \cite{DelSch06}, sect.\,6.6 and \cite{Del99} for expositions), where an equivalent martingale measure is determined by the solution of some optimization problem. 

Our approach is based on the following general statement. Let $X$ be a Banach lattice (\cite{AliBor06}, chap.9) with the topological dual $X^*$. Denote by $X^*_+$ the non-negative cone of $X^*$ and consider a convex cone $C\subset X$. If an element $f\in X^*$ is bounded from above on a certain subset of $C$:
\begin{displaymath}
\sup_{w\in C_1} f(w)<\infty,\ \ \ C_1=\{w\in C:\|w^-\|\le 1\},\ \ w^-=\max\{-w,0\},
\end{displaymath}
then there exists an element $g\in X^*$ such that $g\ge f$, $g(w)\le 0$, $x\in C.$ 
The results of this type were exploited in \cite{Rok05}, \cite{Rok08b} to prove new versions of the Kreps-Yan theorem and in \cite{RokSch06}, \cite{Rok08a} to analyse lower bounds of martingale measure densities. 

Let $S$ be an integrable process, $X=L^1$ and let $K_n\subset L^1$ be the subspace of the elements $(\gamma_n,\Delta S_n)$, where $\gamma_n$ is an $\mathcal F_{n-1}$-measurable vector with bounded components. The topological dual of $L^1$ is the space $L^\infty$ of essentially bounded functions. Assuming the mentioned result, we see that it is sufficient to present a strictly positive element $f\in L^\infty$ with
\begin{displaymath}
\sup\{\mathsf E (wf):\mathsf E( w^-)\le 1,\ w\in K_n\}<\infty
\end{displaymath}
to prove the existence of an equivalent martingale measure with bounded density for $(S_{n-1},S_n)$. Namely, this density corresponds to an element  $g\in L^\infty$: $g\ge f$, $\mathsf E(wg)=0$, $w\in K_n$, indicated above. 

In contrast to \cite{Sch92}, \cite{KabKra94}, \cite{KabStr01}, this approach does not require to prove the closedness (in probability or in $L^1$) of the cone $C$ or the subspace $K_n$. As compared to \cite{DalMorWil90}, \cite{JacShi98}, we do not consider set-valued mappings, whose values lie in infinite dimensional spaces, and do not use delicate Aumann's measurable selection theorem (\cite{AliBor06}, Theorem 18.26) and "projection theorem" (\cite{AliBor06}, Theorem 18.25).   
The mathematical tools of this paper include: (a) relatively simple results of the theory of measurable set-valued mappings with values in $\mathbb R^d$ ("preservation of measurability" theorems, the Kuratowski--Ryll-Nardzewski measurable selection theorem, Castaing representation: sect. 14A, 14B of \cite{RocWet98}); (b) standard duality results (separation theorems, weak$^*$-compactness of the unit ball in the dual space); (c) well-known probabilistic results, concerning the conditional distributions and conditional expectations.

We do not pretend that the proof presented below is simpler, shorter or better than the existing ones. Rather it gives a somewhat different view of the problem.

\section{Auxillary statements}
Let $\mathscr H$ be a sub-$\sigma$-algebra of $\mathscr F$. Recall that a set-valued mapping $F$, assigning some set $F(\omega)\subset\mathbb R^d$ to each $\omega\in\Omega$, is called $\mathscr H$-\emph{measurable}, if $\{\omega: F(\omega)\cap V\neq\varnothing\}\in\mathscr H$ for any open set $V\subset\mathbb R^d$.
A function $\eta:\Omega\mapsto\mathbb R^d$ is called a \emph{selector} of $F$, if $\eta(\omega)\in F(\omega)$ for all  $\omega\in\dom F:=\{\omega':F(\omega')\neq\varnothing\}$. Denote by $\mathscr S(\mathscr H,F)$ the set of $\mathscr H$-measurable selectors of $F$. According to the Kuratowski--Ryll-Nardzewski theorem we have $\mathscr S(\mathscr H,F)\neq\varnothing$ for any $\mathscr H$-measurable set-valued mapping $F$ with closed values $F(\omega)$ (\cite{RocWet98}, Corollary 14.6). Moreover, there exists a \emph{Castaing representation} of $F$, that is, a sequence $(\eta_i)_{i=1}^\infty$ of $\mathscr H$-measurable selectors of $F$ such that the sets $\{\eta_i(\omega)\}_{i=1}^\infty$ are dense in $F(\omega)$ for all $\omega\in\dom F$ (\cite{RocWet98}, Theorem 14.5). 

Let $L^1(\mathscr H,\mathsf P)$ and $L^\infty(\mathscr H,\mathsf P)$ be the Banach spaces of equivalence classes of $\mathscr H$-measurable real-valued functions with the norms $\|w\|_1=\mathsf E|w|$ and $\|w\|_\infty=\esssup|w|$. Denote by $L^p(\mathscr H,\mathsf P,F)$, $p\in\{1,\infty\}$ the set of equivalence classes of $\mathscr H$-measurable vectors, satisfying the conditions $\eta\in F$ a.s., $|\eta|\in L^p(\mathscr H,\mathsf P)$, where $|x|=(x,x)^{1/2}$; and denote by $L^p_+(\mathscr H,\mathsf P)$, $L^p_{++}(\mathscr H,\mathsf P)$, $p\in\{1,\infty\}$ the sets of non-negative and strictly positive elements of $L^p(\mathscr H,\mathsf P)$ respectively. In what follows, for brevity, we omit the argument $\mathsf P$ in the above notation.

By $\ri A$, $\conv A$, $\cone A$ and $\lin A$ we denote the relative interior, the convex hull, the conic hull and the linear span of a subset $A$ of a finite dimensional space. If $A$ is a cone then $A^\circ=\{y\in\mathbb R^d:(x,y)\le 0\}$ is the polar cone. 

The support $\varkappa_\xi$ of the regular conditional distribution $\mathsf P_\xi(\omega,dx)$ of the random vector $\xi$ with respect to $\mathscr H$ is an $\mathscr H$-measurable set-valued mapping:
\begin{displaymath}
\{\omega:\varkappa_\xi(\omega)\cap V\neq\varnothing\}=\{\omega:\mathsf P_\xi(\omega,V)>0\}\in\mathscr H
\end{displaymath}
for any open set $V\subset\mathbb R^d$. The "preservation of measurability" results (\cite{RocWet98}, Proposition 14.11, Example 14.12) show that the set-valued mapping
\begin{equation}\label{eq1}
\omega\mapsto G(\omega)=\lin\varkappa_\xi(\omega)\cap (-\cone \varkappa_\xi(\omega))^\circ\cap \mathbb S_1^d, 
\end{equation}
where $\mathbb S_1^d$ is the unit sphere of $\mathbb R^d$, is $\mathscr H$-measurable. It is easy to see that $G(\omega)\neq\varnothing$, if and only if $\omega$ belongs to the set $A_\xi=\{\omega:0\not\in\ri(\conv\varkappa_\xi(\omega))\}.$ Hence, $A_\xi\in\mathscr H$.
\begin{lemma} \label{lem1}
Let $\xi\in\mathscr S(\mathscr F,\mathbb R^d)$ and $\mathsf P(A_\xi)>0$, where $A_\xi=\{0\not\in\ri(\conv\varkappa_\xi)\}$. Then there exists an $\mathscr H$-measurable random vector $\gamma:\Omega\mapsto\mathbb R^d$ such that $\mathsf P((\gamma,\xi)\ge 0)=1$, $\mathsf P((\gamma,\xi)>0)>0$.
\end{lemma}
\textbf{Proof.} By virtue of the Kuratowski--Ryll-Nardzewski theorem we can take $\eta\in\mathscr S(\mathscr H,G)$, where $G$ is defined by (\ref{eq1}). Put $\gamma=\eta I_{A_\xi}$. Then
\begin{displaymath}
\mathsf P(\gamma,\xi)\ge 0)=\mathsf E\mathsf P_\xi(\omega,\{x:(\gamma(\omega),x)\ge 0\})=1,
\end{displaymath}
since $(\gamma(\omega),x)\ge 0$, $x\in\varkappa_\xi(\omega)$ and $\mathsf P_\xi(\omega,\{x:(\gamma(\omega),x)\ge 0\})=1$ for all $\omega\in\Omega$. Moreover,
\begin{displaymath}
\mathsf P((\gamma,\xi)>0)=\mathsf E\mathsf P_\xi(\omega,\{x:(\gamma(\omega),x)> 0\})>0.
\end{displaymath}
Actually, for any $\omega\in A_\xi$ there exists $y\in\varkappa_\xi(\omega)$ such that $(\gamma(\omega),y)>0$.
Thus, $\mathsf P_\xi(\omega,\{x:(\gamma(\omega),x)>0\})>0$ for all $\omega\in A_\xi$. $\square$

It is clear that Lemma \ref{lem1} leads to the proof of the assertion (i) $\Longrightarrow$ (iii) of Theorem 1. The remaining reasoning of this section prepare the proof of the key implication (iii) $\Longrightarrow$ (ii). In fact, it is sufficient to establish the next result.
\begin{lemma}\label{lem2}
Let $\xi\in L^1(\mathscr F,\mathbb R^d)$ and $0\in\ri(\conv\varkappa_\xi)$ a.s. Then there exists $g\in L^\infty_{++}(\mathscr F)$ such that $\mathsf E(g\xi|\mathscr H)=0$.
\end{lemma}

We prove Lemma \ref{lem2} by the following scheme. Any element $f\in L^1(\mathscr F)$ induces the linear functional on $L^\infty(\mathscr F)$ by the formula $\langle w,f\rangle=\mathsf E(wf)$. Consider the subspace
\begin{equation}\label{eq2}
K=\{(\gamma,\xi):\gamma\in L^\infty(\mathscr H,\lin\varkappa_\xi)\}\subset L^1(\mathscr F).
\end{equation}
Under the assumption $0\in\ri(\conv\varkappa_\xi)$ a.s. there exists an element $f\in L^1_{++}(\mathscr F)$, which is bounded from above one the special subset $K_1$ of $K$:
\begin{displaymath}
 \sup_{w\in K_1} \mathsf E (wf)<\infty,\ \ \ K_1=\{w\in K:\|w^-\|_1\le 1\},\ \ w^-=\max\{-w,0\}
\end{displaymath}
(Lemma \ref{lem5} below). It follows that there exists $g\in L^1(\mathscr F)$, $g\ge f$ such that: $\langle w,g\rangle=0$, $w\in K$ (Lemma \ref{lem7}). This element $g$ satisfies the conditions of Lemma \ref{lem2}.  

We turn to the realization of this scheme. Let  $\xi\in L^1(\mathscr F,\mathbb R^d)$ and $B_\xi=\{\omega:\int |x|\,d\mathsf P_\xi(\omega,dx)<\infty\}$. Clearly, $B_\xi\in\mathscr H$ and $\mathsf P(B_\xi)=1$. We put 
\begin{displaymath}
\psi(\omega,h)=I_{B_\xi}(\omega)\int_{\mathbb R^d} (h,x)^-\,\mathsf P_\xi(\omega,dx) 
\end{displaymath}
and introduce the set-valued mapping
\begin{equation}\label{eq3}
\omega\mapsto T(\omega)=\{h:\psi(\omega,h)\le 1\}\cap (\lin\varkappa_\xi(\omega)).
\end{equation}

\begin{lemma} \label{lem3}
Let $\xi\in L^1(\mathscr F,\mathbb R^d)$ and $0\in\ri(\conv\varkappa_\xi)$ a.s.
Then the set-valued mapping $T$, defined by (\ref{eq3}), is $\mathscr H$-measurable and has compact values $T(\omega)$ a.s.
\end{lemma}
\textbf{Proof.}  The function $\psi:\Omega\times\mathbb R^d\mapsto\mathbb R$ is convex in $h$ for each $\omega$ and $\psi(\omega,0)=0$. Thus the set $T'(\omega)=\{h\in\mathbb R^d:\psi(\omega,h)\le 1\}$ is convex and has non-empty interior for all $\omega\in\Omega$. It follows that the $\mathscr H$-measurability of $T'$ is implied by the simple test of having
$$ \{\omega:x\in T'(\omega)\}=\{\omega:\psi(\omega,x)\le 1\}\in\mathscr H$$
for every $x\in\mathbb R^d$ (\cite{RocWet98}, Example 14.7), and $T$ is $\mathscr H$-measurable as an intersection of $\mathscr H$-measurable set-valued mappings (\cite{RocWet98}, Proposition 14.11(a)).

Assume that $\omega\in A_\xi\cap B_\xi$. Since $0\in\ri(\conv\varkappa_\xi(\omega))$, it follows that for $h\in \lin\varkappa_\xi(\omega)\backslash 0$ the set $\varkappa_\xi(\omega)$ is not contained in the half-space $\{x\in \lin\varkappa_\xi(\omega):(h,x)\ge 0\}$. Thus $\psi(\omega,h)>0$ and the set $T(\omega)$ is compact, because 
$\psi(\omega,h)\to \infty$ as $|h|\to\infty$, $h\in \lin\varkappa_\xi(\omega)$. $\square$

\begin{lemma} \label{lem4}
Let $\xi\in\mathscr S(\mathscr F,\mathbb R^d)$ and $0\in\ri(\conv \varkappa_\xi)$. If $(\gamma,\xi)\ge 0$ a.s. for some $\gamma\in\mathscr S(\mathscr H,\lin\varkappa_\xi)$, then $\gamma=0$ a.s. 
\end{lemma}
\textbf{Proof.} Put $A=\{\gamma\neq 0\}$. For any $\omega\in A$ there exists $y\in \varkappa_\xi(\omega)$ such that $(\gamma(\omega),y)<0$ and hence $\int (\gamma(\omega),x)^-\,\mathsf P_\xi(\omega,dx)>0$. If $\mathsf P(A)>0$, then we obtain the contradiction:
\begin{displaymath}
\mathsf E(\gamma ,\xi)^- = \mathsf E \mathsf E(I_A (\gamma,\xi)^-|\mathscr H)=\mathsf E\left( I_A \int_{\mathbb R^d} (\gamma(\omega),x)^-\,\mathsf P_\xi(\omega,dx)\right)>0. \ \ \square
\end{displaymath}

Let $(\eta_i)_{i=1}^\infty$ be an $\mathscr H$-measurable Castaing representation of $T$ and let $\zeta:\Omega\mapsto\mathbb R^d$ be an $\mathscr H$-measurable vector. Denote by $s(x|A)=\sup\{(x,y):y\in A\}$ the \emph{support function} of a set $A$. From
\begin{displaymath}
 s(\zeta(\omega)|T(\omega))=\sup_{i\in\mathbb N}(\zeta(\omega),\eta_i(\omega))
\end{displaymath}
it follows that the function $s(\zeta|T)$ is $\mathscr H$-measurable. In addition it is a.s. finite, owing to the compactness of $T(\omega)$.
\begin{lemma} \label{lem5}
Let $\xi\in L^1(\mathscr F,\mathbb R^d)$ and $0\in\ri(\conv\varkappa_\xi)$ a.s. For 
\begin{equation} \label{eq4}
f=(1+s(\mathsf E(\xi|\mathscr H)|T))^{-1}\in L^\infty_{++}(\mathscr H)
\end{equation} 
we have
\begin{displaymath}
\beta:=\sup_\gamma\{\mathsf E(\gamma,f\xi): 
 \|(\gamma,\xi)^{-}\|_1\le 1,\ \gamma\in L^\infty(\mathscr H,\lin\varkappa_\xi)\}\le 1.
\end{displaymath}
\end{lemma}
\textbf{Proof.} Put $U^1_+(\mathscr H)=\{g\in L^1_+(\mathscr H): \|g\|_1 \le 1\}$. We see that
\begin{eqnarray*}
 U^1_+(\mathscr F) &=& \{g\in L^1_+(\mathscr F): \mathsf E(\mathsf E(g|\mathscr H )) \le 1\}\\
&=& \bigcup_{w\in U_+^1(\mathscr H) }\{g\in L^1_+(\mathscr F): \mathsf E(g|\mathscr H )\le w\}.
\end{eqnarray*}
Putting $a=\mathsf E(f\xi|\mathscr H)$, we get
\begin{eqnarray*}  
\beta &=& \sup_\gamma\{\mathsf E(\gamma,a): (\gamma,\xi)^{-}\in U^1_+(\mathscr F),\ \gamma\in L^\infty(\mathscr H,\lin\varkappa_\xi)\}\\
&=& \sup_{w\in U_+^1(\mathscr H)} \sup_\gamma\{\mathsf E(\gamma,a): \mathsf E((\gamma,\xi)^{-}|\mathscr H) \le w,\ \gamma\in L^\infty(\mathscr H,\lin\varkappa_\xi)\}.
\end{eqnarray*}

On the set $\{w=0\}$ we have the equality $\mathsf E((\gamma,\xi)^{-}|\mathscr H)=0$. Therefore,
$\mathsf E((\gamma I_{\{w=0\}},\xi)^{-})=0$
and $\gamma I_{\{w=0\}}=0$ by Lemma \ref{lem4}. Putting $\gamma=w \theta$, where $\theta$ is an $\mathscr H$-measurable vector, we obtain
\begin{displaymath}
\beta=\sup_{w\in U_+^1(\mathscr H)} \sup_\theta\{\mathsf E w (\theta,a): \mathsf E((\theta I_{\{w>0\}},\xi)^{-}|\mathscr H)\le 1,\ w \theta\in L^\infty(\mathscr H,\lin\varkappa_\xi)\}.
\end{displaymath}

Since the values of $\theta$ on the set $\{w=0\}$ do not affect $\mathsf E w (\theta,a)$, by the definition of $T$ and the equality $\mathsf E((\theta,\xi)^{-}|\mathscr H)=\psi(\omega,\theta(\omega))$ a.s., we get
\begin{displaymath} 
\beta=\sup_{w\in U_+^1(\mathscr H)} \sup_\theta\{\mathsf E w (\theta,a): \theta\in\mathscr S(\mathscr H,T),\ w\theta\in L^\infty(\mathscr H,\lin\varkappa_\xi)\}.
\end{displaymath}
But $(\theta,a)\le s(a|T)=s(\mathsf E(\xi|\mathscr H)|T)f\le 1$ a.s. for $\theta\in\mathscr S(\mathscr H, T)$. This yields that $\beta\le \sup_{w\in U_+^1(\mathscr H)}\mathsf Ew = 1$. $\square$

Denote by $U^\infty$ the unit ball of the space $L^\infty(\mathscr F)$ and put $U_+^\infty=\{w\in L^\infty_+(\mathscr F):w\in U^\infty\}$, $w^+=\max\{w,0\}$.
\begin{lemma} \label{lem6}
For any element $w\in L^1(\mathscr F)$ we have
$$\|w^+\|_1=\sup\{\langle w,z\rangle:z\in U^\infty_+\}.$$
\end{lemma}
\textbf{Proof.} If $\zeta=I_{\{w\ge 0\}}\in U_+^\infty$, then $\langle w,\zeta\rangle=\|w^+\|_1.$
On the other hand,
\begin{displaymath} 
\langle w,z\rangle\le\langle w^+,z\rangle\le\|w^+\|_1,\ \ z\in U^\infty_+.\ \ \square
\end{displaymath}

Recall that the closure of a convex set $A\subset L^1(\mathscr F)$ in the weak topology $\sigma (L^1,L^\infty)$ coincides with its norm closure in $L^1(\mathscr F)$.  
\begin{lemma} \label{lem7}
Let $K$ be a subspace of $L^1(\mathscr F)$ and $f\in L^\infty(\mathscr F)$. If
\begin{displaymath} 
\sup_{w\in K_1}\langle w,f\rangle<\infty, \ \ K_1=\{w\in K:\|w^-\|_1\le 1\},
\end{displaymath} 
then there exists an element $g\in L^\infty(\mathscr F)$, satisfying the conditions
\begin{displaymath} 
\langle w,g\rangle=0,\ \ w\in K; \ \ g\ge f.
\end{displaymath} 
\end{lemma}
\textbf{Proof.} Put $\lambda=\sup_{w\in K_1}\langle w,f\rangle$. If the assertion of lemma is false, then
\begin{displaymath} 
(f+\lambda U^\infty_+)\cap K^\circ=\varnothing,\ \ \ K^\circ=\{z\in L^\infty(\mathscr F):\langle w,z\rangle\le 0,\ w\in K\}.
\end{displaymath} 
By applying the separation theorem (\cite{AliBor06}, Theorem 5.79) to the $\sigma (L^\infty,L^1)$-compact set $f+\lambda U^\infty_+$ and $\sigma (L^\infty,L^1)$-closed set $K^\circ$, we conclude that there exists $v\in L^1(\mathscr F)$ 
such that
\begin{displaymath} 
\sup_{z\in K^\circ}\langle v,z\rangle<\inf\{\langle v,\zeta\rangle:\zeta\in f+\lambda U^\infty_+\}.
\end{displaymath} 
Since $K$ is a subspace it follows that $\langle v,z\rangle=0$, $z\in K^\circ$ and $v\in \cl_1 K$ by the bipolar theorem (\cite{AliBor06}, Theorem 5.103), where $\cl_1 K$ is the closure of $K$ in the norm topology of $L^1(\mathscr F)$. Moreover,
\begin{eqnarray} \label{eq5}
\langle v,f\rangle+\lambda\inf\{\langle v,\eta\rangle:\eta\in U^\infty_+\}>0.
\end{eqnarray}

By Lemma \ref{lem6} we have
\begin{equation} \label{eq6}
\inf\{\langle v,\eta\rangle:\eta\in U^\infty_+\}=-\sup\{\langle -v,\eta\rangle:\eta\in U^\infty_+\}=-\|v^-\|_1.
\end{equation}
If $v^-=0$ then $\langle v,f\rangle>0$ and $\alpha v\in L^1_+\cap\cl_1 K$ for any $\alpha>0$. Hence, the functional $w\mapsto \langle w,f\rangle$ is unbounded from above on the ray $\{\alpha v:\alpha>0\}$, which is contained in the set
\begin{displaymath} 
\cl_1 K_1\supset\cl_1\biggr(\{w:\|w^-\|_1<1\}\cap K\biggl)\supset \{w:\|w^-\|_1<1\}\cap \cl_1 K.
\end{displaymath} 
Here we have used an elementary inclusion $\cl_1(A\cap B)\supset A\cap\cl_1 B$, which holds true when the set $A$ is open.

Thus, $\|v^-\|_1>0$ and it follows from (\ref{eq5}), (\ref{eq6}) that 
\begin{displaymath} 
\langle v/\|v^-\|_1,f\rangle>\lambda
\end{displaymath} 
in contradiction to the definition of $\lambda$, since $v/\|v^-\|_1\in K_1$. $\square$

Note that for any $f\in L^1(\mathscr F)$ there exists an $\mathscr H\otimes\mathscr B(\mathbb R^d)$-measurable function $\varphi(\omega,x)$, satisfying the condition $\mathsf E(f|\mathscr H\vee\sigma(\xi))=\varphi(\omega,\xi)$ a.s. Here $\mathscr B(\mathbb R^d)$ is the Borel $\sigma$-algebra of $\mathbb R^d$ and $\sigma(\xi)$ is the $\sigma$-algebra, generated by $\xi$. 
This statement follows from the fact that the $\sigma$-algebra $\mathscr H\vee\sigma(\xi)$ is generated by the mapping $\omega\mapsto(\omega,\xi(\omega))$ from $\Omega$ to the measurable space $(\Omega \otimes\mathbb R^d,\mathscr H\otimes\mathscr B(\mathbb R^d))$ (see \cite{AliBor06}, Theorem 4.41).  

\textbf{Proof of Lemma \ref{lem2}.} Let $K$ be the subspace, defined by (\ref{eq2}). By Lemma \ref{lem5} the element  $f$, defined by (\ref{eq4}), considered as a functional on $L^1(\mathscr F)$, satisfies the conditions of Lemma \ref{lem7}. An element $g$, indicated in Lemma \ref{lem7}, belongs to $L^\infty_{++}(\mathscr F)$ and
\begin{equation} \label{eq7}
\langle (\gamma,\xi),g\rangle=\mathsf E(\gamma,\mathsf E(g\xi|\mathscr H))=0,\ \ \gamma\in L^\infty(\mathscr H,\lin\varkappa_\xi).
\end{equation}
But (\ref{eq7}) is reduced to the equality $\mathsf E(g\xi|\mathscr H)=0$. Indeed, since 
\begin{displaymath} 
\mathsf E(g\xi|\mathscr H)=\mathsf E(\mathsf E(g|\mathscr H\vee\sigma(\xi))\xi|\mathscr H)=\int\varphi(\omega,x)x\,\mathsf P_\xi(\omega,dx)\in \lin\varkappa_\xi\ \textrm{a.s.},
\end{displaymath} 
where $\varphi$ is some $\mathscr H\otimes\mathscr B(\mathbb R^d)$-measurable function, by putting 
\begin{displaymath} 
\gamma=\mathsf E(g\xi|\mathscr H) I_{\{|\mathsf E(g\xi|\mathscr H)|\le M\}}\in L^\infty(\mathscr H,\lin\varkappa_\xi)
\end{displaymath} 
and passing in (\ref{eq7}) to the limit as $M\to\infty$, we conclude that $\mathsf E(g\xi|\mathscr H)=0$ by the monotone convergence theorem. $\square$

\section{Proof of the Dalang-Morton-Willinger theorem} 
(iii) $\Longrightarrow$ (ii). Let us first show that the supports $\varkappa_n$ of the regular conditional distributions $\mathsf P_n(\omega,dx)$ are a.s. invariant under equivalent changes of measure. Let $\mathsf P'$ be an equivalent (to $\mathsf P$) probability measure and let $(Z_n')_{n=0}^N$ be the correspondent density process:
\begin{displaymath} 
Z_n'=\mathsf E\left(\frac{d\mathsf P'}{d\mathsf P}\biggr|\mathscr F_n\right),\ \ n=0,\dots,N.
\end{displaymath} 
Denote by $\mathsf P_{n-1}'(\omega,dx)$ the regular conditional distribution of $\Delta S_n$ with respect to $\mathcal F_{n-1}$, induced by the measure $\mathsf P'$. By the Bayes formula (\cite{Shi98}, chap. V, \S 3a) we have
\begin{displaymath} 
\mathsf P_{n-1}'(\omega,B)=\mathsf E_{\mathsf P'}(I_{\{\Delta S_n\in B\}}|\mathscr F_{n-1})=\frac{\mathsf E(I_{\{\Delta S_n\in B\}} Z_n'|\mathscr F_{n-1})}{Z_{n-1}'}\ \ \textrm{a.s.},
\end{displaymath} 
where $B$ is a Borel subset of $\mathbb R^d$. Let $\varphi_n$ be an $\mathscr F_{n-1}\otimes\mathscr B(\mathbb R^d)$-measurable function, satisfying the condition
\begin{displaymath} 
\varphi_n(\omega,\Delta S_n(\omega))=\mathsf E(Z_n'|\mathscr F_{n-1}\vee\sigma(\Delta S_n))(\omega)\ \textrm{a.s.}
\end{displaymath} 
Then
\begin{displaymath} 
\mathsf P_{n-1}'(\omega,B)=\frac{\mathsf E(I_{\{\Delta S_n\in B\}} \varphi_n(\omega,\Delta S_n)|\mathscr F_{n-1})}{Z_{n-1}'}=\frac{\int_B\varphi_n(\omega,x)\,\mathsf P_{n-1}(\omega,dx)}{Z_{n-1}'(\omega)}\ \textrm{a.s.}
\end{displaymath} 
For any $A_{n-1}\in\mathcal F_{n-1}$ we have
\begin{displaymath} 
\mathsf E(\varphi_n^-(\omega,\Delta S_n)I_{A_{n-1}})=\mathsf E(I_{A_{n-1}}(\omega)\int_{\mathbb R^d} \varphi_n^-(\omega,x)\,\mathsf P_{n-1}(\omega,dx))=0,
\end{displaymath} 
since $\varphi_n(\omega,\Delta S_n(\omega))>0$ a.s. Hence, for some set $\Omega'$ with $\mathsf P(\Omega')=1$ the inequality $\varphi(\omega,x)>0$ holds true almost everywhere with respect to the measure $\mathsf P_{n-1}(\omega,dx)$ if $\omega\in\Omega'$. This means that for $\omega\in\Omega'$ the measures $\mathsf P_{n-1}(\omega,dx)$, $\mathsf P_{n-1}'(\omega,dx)$ are equivalent and their supports are the same. 

By virtue of the proved invariance property we can assume that $S_n\in L^1(\mathscr F_n,\mathsf P)$. Otherwise we may pass to the measure $\mathsf P'$ with density 
\begin{displaymath}
d\mathsf P'/d\mathsf P=c \exp(-\sum_{n=0}^N |S_n|),
\end{displaymath}
 where $c>0$ is the normalizing constant. In this case $S_n\in L^1(\mathscr F_n,\mathsf P')$ and $d\mathsf Q/d\mathsf P=d\mathsf Q/d\mathsf P'\cdot d\mathsf P'/d\mathsf P\in L^\infty(\mathscr F)$, if $d\mathsf Q/d\mathsf P'\in L^\infty(\mathscr F)$.

By Lemma \ref{lem2} for any $n=1,\dots,N$ there exists $g_n\in L^\infty_{++}$ such that 
\begin{displaymath} 
\mathsf E(g_n \Delta S_n|\mathscr F_{n-1})=0.
\end{displaymath} 
Further argumentation is borrowed from \cite{DelSch06}, sect.\,6.7.
If $N=1$, then the measure $\mathsf Q$ with density $d\mathsf Q/d\mathsf P=g_N/\mathsf E g_N$ has the desired property:
\begin{displaymath} 
\mathsf E_\mathsf Q(I_A \Delta S_N)=\mathsf E(I_A g_N\Delta S_N)/\mathsf E g_N=0,\ \ A\in\mathscr F_{N-1}.
\end{displaymath} 

Assume now that the assertion under consideration is true for $N=m-1$. Applying it to the processes $(S_n)_{n=1}^m$, $(S_n)_{n=0}^1$, we see that on $\mathscr F_m$ there exists an equivalent martingale measure $\mathsf Q'$ for $(S_n)_{n=1}^m$ with $\mathsf P$-density  $d\mathsf Q'/d\mathsf P\in L^\infty(\mathscr F_m)$ and on $\mathscr F_1$ there exists an equivalent martingale measure $\mathsf Q_1$ for $(S_n)_{n=0}^1$ with $\mathsf Q'$-density $f_1=d\mathsf Q_1/d\mathsf Q'\in L^\infty(\mathscr F_1)$. Define on $\mathscr F_m$ the probability measure $\mathsf Q$ by
\begin{displaymath} 
\mathsf Q(A)=\int_A f_1\,d\mathsf Q',\ \ A\in\mathscr F_m. 
\end{displaymath} 
Clearly, $d\mathsf Q/d\mathsf P=f_1 d\mathsf Q'/d\mathsf P\in L^\infty_{++}(\mathscr F_m)$. It remains to check that the process $(S_n)_{n=0}^m$ is a $\mathsf Q$-martingale. 

For $n=1$ this follows from the definition of $f_1$:
\begin{displaymath} 
\mathsf E_\mathsf Q(I_A \Delta S_1)=\mathsf E_{\mathsf Q'}(I_A f_1\Delta S_1)=\mathsf E_{\mathsf Q_1}(I_A \Delta S_1)=0, \ \ A\in\mathscr F_0.
\end{displaymath} 
For $n>1$ by the definition of $\mathsf Q'$ and the $\mathscr F_1$-measurability of $f_1$ we have:
\begin{displaymath} 
\mathsf E_\mathsf Q(I_A \Delta S_n)=\mathsf E_{\mathsf Q'}(I_A f_1\Delta S_n)=\mathsf E_{\mathsf Q'}(I_A f_1\mathsf E_{\mathsf Q'}(\Delta S_n|\mathscr F_{n-1}))=0,\ \ A\in\mathscr F_{n-1}.
\end{displaymath} 

(i) $\Longrightarrow$ (iii). Suppose that $0\not\in\ri(\conv\varkappa_{n-1}(\omega))$ for some $n$ on a set $A\in\mathscr F_{n-1}$ with $\mathsf P(A)>0$. Put $\gamma_j=0$, $j\neq n$. Applying Lemma \ref{lem1} to $\xi=\Delta S_n$ we conclude that there exists an $\mathscr F_{n-1}$-measurable vector $\gamma_n:\Omega\mapsto\mathbb R^d$ such that
\begin{displaymath} 
\mathsf P(G_N^\gamma=(\gamma_n,\Delta S_n)\ge 0)=1,\ \ \ \mathsf P(G_N^\gamma=(\gamma_n,\Delta S_n)>0)>0.
\end{displaymath} 
Thus the NA condition is violated.

(ii) $\Longrightarrow$ (i). If $S$ fails the NA property, then for some $n$ there exists $\gamma_n\in L^\infty(\mathscr F_{n-1},\mathbb R^d)$ such that
\begin{equation}\label{eq8}
\mathsf P((\gamma_n,\Delta S_n)\ge 0)=1, \ \ \ \mathsf P((\gamma_n,\Delta S_n)>0)>0,
\end{equation}
that is, there exists a one-step arbitrage opportunity. We do not reproduce here the simple and well-known proof of this statement (see, e.g., \cite{Sch92}, Lemma 1.2). Evidently, (\ref{eq8}) contradicts the existence of an equivalent martingale measure $\mathsf Q$, since $\mathsf E_\mathsf Q(\gamma_n,\Delta S_n)=\mathsf E_\mathsf Q(\gamma_n,\mathsf E_\mathsf Q(\Delta S_n|\mathscr F_{n-1}))=0$. $\square$


\end{document}